\documentclass[11pt,a4paper]{amsart}
\usepackage{graphicx}

\usepackage{amsmath}%
\usepackage{amssymb,indentfirst,calc,mathrsfs,euscript,bbm}%

\usepackage{setspace}%
\usepackage[reals]{layout}%
\usepackage{xr}%
\usepackage{amscd}%
\usepackage{comment}
\usepackage{epstopdf}

\def\square{\hfill\hbox{\vrule height .9ex width .8ex depth -.1ex}}

\def\rit{\mathbb R}

\def\val{\hbox{\tt val}}

 \def\P{{\mathsf P}} %Probability
 \def\E{{ \mathsf E}} % Esperance

\def\val{\hbox{\tt val}}

\def\U{{\bf U}}

\def\V{{\bf V}}

\def\bx{{\bf x}}

\def\bk{{\bf k}}
\def\bg{{\bf g}}
\def\bq{{\bf q}}
\def\bx{{\bf x}}
\def\by{{\bf y}}
 \def\V{{\mathsf V}} 
  \def\G{{\mathsf G}} 
   \def\v{{\mathsf v}} 
    \def\W{{\mathsf W}} 
         
            \def\U{{\mathsf U}}

\newtheorem{defi}{Definition}[section]

\newtheorem{cor}{Corollary}[section]
\newtheorem{pro}{Proposition}[section]

\begin{document}

\title[ROS]{
Limit value of dynamic   zero-sum games with vanishing stage duration  }
\author{
Sylvain Sorin }
\address{
Sorbonne UniversitŽ\'es, UPMC Univ Paris 06, Institut de MathŽ\'ematiques de Jussieu-Paris Rive Gauche, UMR 7586, CNRS, Univ Paris Diderot, Sorbonne Paris Cit\'eŽ, F-75005, Paris, France }

\email{sylvain.sorin@imj-prg.fr}%
\date{ March  2016. Dedicated to the memory of L.S. Shapley. A fisrt version of this work was presented at the conference "Mathematical Aspects of Game Theory and Applications",
Roscoff, June 30 - July 4, 2014. This research was partially supported by  PGMO  2014-LMG}

\bibliographystyle{apalike}

\begin{abstract}
We consider  two person zero-sum games where the players control, at discrete times $\{t_n \}$ induced by  a partition $ \Pi $ of $\rit^+ $,  a continuous time  Markov  state process. 
We prove that the limit of the values $v_{\Pi}$ exist as the mesh of $\Pi$ goes to 0.
The analysis covers  the cases of :\\
 1)  stochastic games (where both players know the state)\\
2) symmetric no  information.\\
The proof is by reduction to a deterministic differential game.

\end{abstract}

\maketitle

\section{Introduction}
Repeated interactions in a stationary environment have been traditionally represented by dynamic  games played in stages. 
An alternative approach is to consider  a continuous time process on which the players act at discrete times. 
In the first case the expected number of  interactions increases as  the weight  $\theta_n$ of each stage $n$  goes to zero.
In the second case  the number of interactions increases when the duration  $ \delta_n$ of each time interval $n$ vanishes. 

In a repeated game framework  one can  normalize the model  using  the evaluation   of the stages, so that stage $n$ is associated to time $t_n= \sum_{j=1}^{n-1}\theta_j$, and then consider the game played on $[0,1]$ where time $t$ corresponds to the fraction $t$ of the total duration. Each evaluation $\theta$ (in the original repeated game) thus induces  a partition  $\Pi_\theta$ of $[0,1]$ with vanishing mesh  corresponding to vanishing stage weight.
Tools adapted from continuous time models can  be used to obtain convergence  results, given an ordered set   of evaluations,  for  the corresponding family of values $v_\theta$, see e.g. for different classes of games,  Sorin \cite{Sbm1}, \cite{SFC}, \cite{SDGA}, Vieille \cite{V}, Laraki \cite{LR2}, Cardaliaguet, Laraki and Sorin \cite{CLS}.

In the alternative approach considered here, there is a given evaluation $\bk$ on $\rit^+$ and one consider  a sequence  of partitions $ \Pi(m) $ of $ \rit ^+ $ with vanishing mesh corresponding to vanishing stage duration and the associated sequence of values.

In both cases, for each given partition the value function exists  at the times defined by the partition and  the stationarity of the model allows to write a recursive formula ($RF$). Then one  extends the value function to $[0,1]$ (resp. $\rit^+$)  by linearity and  one considers   the family of values  as the mesh of the partition goes to 0. 
The next two steps in the proof of convergence of the family of values consist in defining  a PDE  (Main Equation $ME$) and   proving :\\
1) that any accumulation point  of the family is a viscosity solution  of ($ME$) \\
2) that ($ME$) has a unique viscosity solution.

Altogether  the tools are quite similar to those used in differential games, however in the current framework the state is basically a random variable and the players use mixed strategies.

Section 2 describes the model.  Section 3 is devoted to the framework  where both players observe the state variable. Section 4 deals with the situation where the state is unknown but the moves are observed. In both cases the analysis is done by reduction to a differential game.  Section 5 presents the main results concerning differential games that are used in the paper. 
\section{Smooth continuous time games and discretization}

\subsection{Discretization of 	a continuous time process and associated game }$ $

Consider  a time homogeneous state process  $Z_t$,  defined on $\rit ^+ =[0, + \infty)$,  with values in a state space $\Omega$ and 
an evaluation given by a   probability  density $\bk(t) $ on $\rit^+$.\\
 Each  partition  $\Pi =\{ t_1=0, t_2, ..., t_n, ...\} $ of $\rit^+$ induces a discrete time game as follows. 
 The time interval $L_n= [t_{n}, t_{n+1}[$   corresponds to stage $n$ and   has duration  $\delta_n $.
 %  and the sequence  $\{\delta_n\}$ is assumed to be decreasing {\it  why??}. \\
The law of $Z_t$ on $L_n$   is determined by its value at time $t_n$, $\hat Z_n= Z_{t_{n}}$ and the actions 
%$\hat x_n$, $\hat y_n$  ,
$(i_n, j_n) \in I \times J $ chosen by the players at time $t_{n}$,  that last for   stage $n$.
% (hence $(i_t, j_t) = (i_n, j_n) $  for $t\in L_n$).  \\
The payoff at time $t$ in stage $n$ ($t \in L_n$)  is  defined  trough a map $ \bg$ from $\Omega \times I\times J $ to $\rit$:
$$
{ \bf g} _\Pi (t) = { g}(Z_t,  i_n, j_n)
 $$
{\it (An  alternative choice leading to the same  asymptotic results would be  $ g _\Pi (t) = { \bf g}( \hat Z_n,  i_n, j_n) $).}\\
The  evaluation along a play is: 
$$
\gamma_{\Pi, \bk} =  \int_0^{+\infty} {\bf g}_\Pi (t) \bk(dt)
$$
and  the  corresponding value  function %, defined on $\Omega$ 
is  $v_{\Pi, \bk}$.\\ 
One will study  the asymptotics  of the family  $\{v_{\Pi, \bk}\}$ as the mesh $\delta =\sup  \delta_n $
of the partition $\Pi$  vanishes.

\subsection{Markov process }$ $

From now on  we consider the case where $Z_t, t  \in \rit ^+$ follows a continuous time Markov process: it is specified by a   transition rate $\bf q$ that belongs to the  set $ {\mathcal M}$ of  real   bounded maps   on  $ I\times J  \times\Omega \times\Omega $    with: 
$$
 { \bf q} (  i, j) [ \omega, \omega '] \geq 0, \quad  \mbox{ if $ \omega' \not= \omega$,\quad and } 
\sum_{\omega ' \in \Omega} {\bf q}( i, j) [ \omega,  \omega '] = 0, \quad \forall i, j,  \omega.
$$
%For $ ({\bf i, \bf j} ) \in I^{\Omega}\times  J^{\Omega}$ define 
%$$
%{\bf q}({\bf i, \bf j} )  [ \omega, \omega ']  = q_{{\bf i }(\omega), {\bf j}(\omega ')}[\omega, \omega ']
%$$
Let $\P^h(i, j ), h\in \rit^+ $ be the continuous  time Markov chain on $ \Omega$  generated by the kernel ${\bf q}({ i,  j} ) $:
$$
 \dot \P ^h ({ i,  j} ) = \P ^h (i, j )\,  {\bf q}(i, j )  = {\bf q}(i, j ) \; \P ^h(i, j )
 $$
and for $t\geq 0$ :
$$
\P^{t + h}(i, j  )  =  \P^{t}(i, j  )\, e^{ h \; {\bf q} (i, j )  }.
$$
In particular, one has:
\begin{eqnarray*}
\P^{h}(i, j)   [z, z'] &=& Prob \, ( Z_{t+h}= z'| Z_t = z
%,{ i}_s= {i},{  j}_s= {j }, 0\leq t \leq s\leq t + h 
), \qquad \forall t\geq 0 \cr
&=& {\bf 1}_{\{z\}} (z')  +  h \,  {\bf q} (i, j  )  [ z,  z '] + o(h)
\end{eqnarray*}

\bigskip

\subsection{Hypotheses }$ $

  One assume from now on:\\
 the state space  $\Omega$ is finite,\\ 
 the evaluation $\bk$ is Lipschitz continuous on $\rit^+$.\\
  the action sets $I$, $J$  are compact metric spaces, \\
  the payoff  ${\bf g}$    and the transition ${\bf q} $ are  
  %bounded.\\ Moreover ${\bf g}$  and ${\bf q} $  are measurable and 
  continuous on $I \times J$.
  
\subsection{Notations }$ $

If $\mu$ is   a  bounded measurable function defined on $I\times J$ with values in a convex set,  $\mu (x,y)$  denotes its multilinear extension to $X\times Y$, with $X = \Delta(I)$ (resp. $Y = \Delta(J)$), set of regular  Borel probabilities  on $I$ (resp. $J$). (This applies in particular to  $\bg$ and $\bq$).\\
For $\zeta \in \Delta (\Omega)$  and   $\mu \in \rit^{\Omega^2} $ we define :
$$
 \zeta *\mu  \;  (z) =  \sum_{\omega \in \Omega}  \zeta (\omega) {\mu}  [ \omega , z] . 
 $$
(When $g$ is a map from $\Omega$ to itself and $\mu [\omega, z]= \bold{1} _{\{{g(\omega) = z}\}}$, $\zeta * g$ is the usual image measure). \\
In particular, if $\zeta_t\in \Delta(\Omega) $ is the law of $Z_t$ one has, if $(i,j)$ is played on $[t, t+h]$
$$
 \zeta_{t+h}= \zeta_t * \P^h (i,j) $$
 and 
 $$
 \dot \zeta_t = \zeta_t  * {\bf q}(i, j).
 $$
 Similarly we use   the following notation for a transition probability or a transition rate  $\mu$  operating on a real  function $f$ on $\Omega$:  
$$
 \mu [z, .]\circ f(\cdot) = \sum_{z'}   \mu[z, z'] f( z') =  \mu \circ f  \, [z]. 
  $$

\section{State controlled and publicly observed}

This section is devoted to the case were the process  $Z_t$ is controlled by both players and observed by both (there is no assumptions on  the signals on the actions). A stage $n$ (time $t_n$) both players know $Z_{t_n}$.
This corresponds to a stochastic game $G$  in continuous time analyzed trough a time discretization along $\Pi$, $G_\Pi$.\\
Previous related papers to stochastic games in continuous time include Zachrisson \cite{Za64},  Tanaka and Wakuta \cite{TW}, Guo and Hernandez-Lerma  \cite{GH1}, \cite{GH2}, Neyman  \cite{N2012}.\\
The approach via time discretization is related to similar procedures in differential games, see Section 5,  Fleming \cite{F57},  \cite{F61},  \cite{F64}, Scarf \cite{S57} and Neyman \cite{N2013}.

%%%%%%%%%%%%%%%%%%%%%%%%%%%%%%%%%%%%%%%%
 
\subsection{General case }$ $ \\
Consider  a general evaluation $\bk$. Since $\bk$ is fixed during the analysis we will write $v_{\Pi}$ for $v_{\Pi, \bk}$, defined on $\rit ^+\times\Omega .$

\subsubsection{Recursive formula}$ $ \\
 The hypothesis on the data implies that $v_\Pi$ exists, see e.g. \cite{MSZ}, Chapters IV and VII,  or \cite{NSSto}, and in the current  framework the recursive formula takes the following form:
\begin{pro} $ $ \\
The  game $G_\Pi$   has a value  $v_{\Pi }$ satisfying the recursive equation:
 \begin{eqnarray}\label{RF1} 
 \nonumber
v_{\Pi }(t_{n}, Z_{t_{n}}) &=& \val_{X\times Y}  \, \E_{ z, x, y}  [\int_{t_{n}}^{t_{n+1} }{\bf g}(Z_s,  i, j)  \bk(s) ds  + 
 v_{\Pi} (t_{n+1}, Z_{t_{n+1}} )]\\
 &=& \val_{X\times Y}  \, [ \E_{ z, x, y}( \int_{t_{n}}^{t_{n+1} }{\bf g}( Z_s,  i, j) \bk(s) ds )  +    \P^{\delta_n} ( x, y)[Z_{t_{n}},.] \circ  v_{\Pi}(t_{n+1}, .)]
\end{eqnarray}
\end{pro}
\noindent\underline{Proof}\\
This is the basic recursive formula for the stochastic game with state space $\Omega$, action sets $I$ and $J$ and  transition kernel $  \P^{\delta_n} ( i, j)$, going back to Shapley \cite{Sh}.
\square\\

Recall that the value  $v_{\Pi }(.,  z)$  is defined at times $t_n \in \Pi$ and extended by  linearity to $\rit^+$.

\subsubsection{Main equation}$ $ \\
The first property is standard in this framework.
\begin{pro}$ $ \\
The family of values  $\{v_{\Pi } \}$  is uniformly Lipschitz w.r.t.  $t \in \rit^+$.  
\end{pro}
Denote thus by ${\bf V}$ the (non empty) set of accumulation points of the family $\{v_{\Pi } \}_{\Pi}$ (for the uniform convergence on compact subsets of $\rit^+ \times \Omega$) as the mesh $\delta$ vanishes.

\begin{defi}
A  continuous real function $u$ on  $\rit ^+\times\Omega $ is a viscosity solution of:
\begin{equation}\label{ME1}
0 = \frac {d}{dt} u( t, z)  + \val_{X\times Y} \{{\bf g}(z, x, y) \bk(t) + {\bf q} ( x, y)[z, .] \circ u (t, \cdot)\},
\end{equation}
if for any  real function $\psi$, ${\mathcal C}^1$ on $\rit ^+\times\Omega $ with $u - \psi$ having a strict maximum at $(\bar t,  \bar z)\in \rit ^+\times\Omega $:
 \begin{equation*}
0 \leq  \frac {d}{dt} \psi( \bar t, \bar z)  + \val_{X\times Y} \{{\bf g}(\bar z, x, y) \bk(\bar t) + {\bf q} ( x, y)[\bar z, .] \circ \psi (\bar t, \cdot)\}
\end{equation*}
and the dual condition.
\end{defi}

\begin{pro}$ $ \\
Any  $u \in {\bf V}$    is a viscosity solution of 
(\ref{ME1}).
\end{pro}
\noindent\underline{Proof}\\
Let $\psi (t,z)$ be a  ${\mathcal C}^1$ test function such that  $ u - \psi$ has a strict maximum at $(\bar t,  \bar z)$. Consider a sequence $V_m = v_{ \Pi(m)}$ converging uniformy  locally to $u$  as $m \rightarrow \infty$ and  let $( t^*(m), z(m)) $ be a minimizing sequence for $  \{( \psi - V_m ) (t, z), t \in \Pi_m \} $. In particular  $( t^*(m), z(m)) $ converges to $(\bar t,  \bar z)$ as $m \rightarrow \infty.$  
Given $x^*_m$ optimal   for $V_m ( t^*(m), z(m))$  in (\ref{RF1}), one obtains,  with  $t^*(m) = t_{n} \in \Pi_m$: 
$$ 
V_m(t_{n}, z(m)) 
\leq 
\E_{z(m),x^*_m,y}  [ \int_{t_{n-1}}^{t_{n} }{\bf g}( Z_s,  i, j) \bk(s) ds  ] +    \P^{\delta_n} ( x^*_m, y)[z(m),.] \circ  V_m (t_{n+1}, .), \qquad \forall y \in Y,
 $$
 so that by the choice of $( t^*(m), z(m))$:
\begin{eqnarray*}
\psi (t_{n}, z(m)) 
&\leq &
\E_{z(m),x^*_m,y}  [ \int_{t_{n-1}}^{t_{n} }{\bf g}( Z_s,  i, j) \bk(s) ds  ]   +    \P^{\delta_n} ( x^*_m, y)[z(m),.]  \circ  \psi (t_{n+1}, .)]\\
&\leq& \delta_n \bk(t_{n}) \; {\bf g} (z(m), x^*_m, y) + \psi (t_{n+1}, z(m) ) + \delta_n \;   {\bf q} ( x^*_m, y) [z(m),.] \circ  \psi (t_{n+1}, .)  + o(\delta_n).
\end{eqnarray*}
This implies:
\begin{eqnarray*}
0 &\leq & \delta_n \bk(t_{n}) \; {\bf g} (z(m), x^*_m, y) + \delta_n \;  \frac{d}{dt} \psi (t_{n}, z(m) ) + \delta_n \;   {\bf q} ( x^*_m, y)[z(m), .]  \circ  \psi (t_{n+1}, .)  + o(\delta_n)
\end{eqnarray*}
hence dividing  by $\delta_n$ and taking the limit  as $m \rightarrow \infty$ one obtains, for some accumulation point $x^* $ in the compact set  $\Delta (I)$: 
\begin{equation}
0 \leq   \bk(\bar t) \; {\bf g} (\bar z, x^*, y) +  \frac{d}{dt} \psi (\bar t, \bar z ) +  {\bf q} ( x^*, y) [\bar z ,.] \circ  \psi (\bar t, .),   \qquad \forall y \in Y 
\end{equation}
so that:
$$
0 \leq    \frac{d}{dt} \psi (\bar t, \bar z ) +  \val_{X\times Y}[  \; {\bf g} (\bar z, x, y)\,   \bk(\bar t)+  {\bf q} ( x, y) [\bar z ,.] \circ  \psi (\bar t, .)].
$$
%hence the result.
 
\square

 %%%%%%%%%%%%%%%%%%%%%%%%%%%%%%%%%%%%%%%%%%%%%%%%%%%%%%%%%%%% 
\subsubsection{Convergence}$ $ \\
A first proof of the convergence of the family   $\{v_{\Pi } \}_{\Pi}$  would follow from the property: \\
$(P)$
Equation (\ref{ME1}) has a unique viscosity solution.\\

%\begin{pro}$ $ \\ 
% Equation (\ref{ME1}) has a unique viscosity solution hence the family of values converge.
%\end{pro}
%\noindent\underline{Proof}\\
%(\ref{ME1}) is an HJB equation with Hamiltonian: 
%$$
%H(t, z, v) =  \val_{X\times Y} [ {\bf g}(z, x, y) \bk(t) + {\bf q} ( x, y)[z, .] \circ v (t, \cdot)].
%$$
%Uniqueness is quite standard in this finite  framework, see Bardi and Capuzzo Dolcetta \cite{BC}, Cardaliaguet \cite{Cii10}.
%\square\\
%

An alternative approach is to relate the game to a differential game on an extended state space $\Delta(\Omega)$.
Define  $ {V} _\Pi$ on $ \rit^+ \times \Delta(\Omega)$ as the expectation of $v_\Pi$, namely: 
$$
 {V} _\Pi(t , \zeta)  = \langle \zeta, v_\Pi(t, .) \rangle = \sum_{\omega\in \Omega} \zeta({\omega})  v_\Pi(t, {\omega})
$$  
 and denote   $ {\bf X} =X^{\Omega}$ and ${\bf Y } = Y^{\Omega}$. 

\begin{pro}$ $ \\
$ {V}_\Pi $ satisfies:
\begin{eqnarray}\label{RF1'}
 V _{\Pi }(t_{n}, \zeta_{t_{n}}) &=& \val_{  {\bf X} \times{\bf Y }}  \,  [ \sum_\omega  \zeta_{t_n}(\omega) \E_{\omega, \bx(\omega), \by(\omega}) (\int_{t_{n}}^{t_{n+1} }  {\bf g}( Z_s, i, j)\rangle  \bk(s)  ds )   + 
V_{\Pi }(t_{n+1}, \zeta_{t_{n+1}})]
\end{eqnarray}
where $\zeta_{t_{n+1}}(z) = \sum _\omega \zeta_{t_{n}}(\omega) \P^{\delta_n} (\bx (\omega), \by(\omega)) (\omega, z)$.
\end{pro}
\noindent\underline{Proof}\\
(\ref{RF1'}) follows from (\ref{RF1}), the definition of $V_\Pi$  and the formula expressing  $\zeta_{t_{n+1}}$. By independence the optimization in  $X$ at each $\omega$ can be replaced by optimization in $\bf X$ and one uses the linearity in the transition.
\square\\

Equation (\ref{RF1}) corresponds to the usual approach following the trajectory of the process. Equation (\ref{RF1'}) expresses the dynamics of the law   $\zeta$ of the process,  where the players act differently at different states $\omega$.\\

%
%Let  $ {V} (t , \zeta)  = \langle \zeta, v(t, .) \rangle = \sum_{\omega} \zeta({\omega})  v(t, {\omega})$  where $v \in \bf V$.
%
%\begin{pro}
%$ {V} (t , \zeta)$ is a viscosity solution of:
%\begin{equation}\label{ME1'}
%0 =  \frac{d}{dt} V (t , \zeta)  + \val_{  {\bf X} \times{\bf Y }} \; [\langle \zeta, {\bf g}( ., {\bf x}(.),   {\bf y}(.))\rangle  \bk(t) +  \langle  f (\zeta,  {\bf x}, {\bf y}), \nabla  V(t , \zeta)   \rangle ]
%\end{equation}
%\end{pro}
%\noindent\underline{Proof}\\
%(\ref{ME1}) is an HJB equation with Hamiltonian: 
%$$
%H(t, z, v) =  \val_{X\times Y} [ {\bf g}(z, x, y) \bk(t) + {\bf q} ( x, y)[z, .] \circ v (t, \cdot)].
%$$
%Uniqueness is quite standard in this framework, see Bardi and Capuzzo Dolcetta \cite{BC}, Cardaliaguet \cite{Cii10}.
%\square\\

\subsubsection{Related differential game}$ $ \\
We will prove that the recursive  equation (\ref{RF1'})
  is satisfied by the value  of  the  time discretization along $\Pi$ of  the mixed extension of a deterministic differential game  $\mathcal G$ (see Section 5)   on $\rit ^+$,  defined as  follows:\\
 1)   the state space is $\Delta(\Omega)$,\\
2)   the action sets are $ {\bf I} =I^{\Omega}$ and ${\bf J } = J^{\Omega}$,\\
 3)  the dynamics on $\Delta(\Omega) \times \rit ^+$  is:
  $$
 \dot \zeta_t  = f (\zeta_t,  {\bf i}, {\bf j})
 $$
 with 
  $$
f (\zeta,  {\bf i}, {\bf j}) (z) = \sum_{\omega \in \Omega}  \zeta_t(\omega)  {\bf q}( {\bf i}(\omega), {\bf j}(\omega)) [ \omega, z].
 $$
4)  the  flow payoff function is given by:
 $$
 \langle \zeta , {\bf g}( ., {\bf i}(.),   {\bf j}(.))\rangle  =   \sum_{\omega \in \Omega}  \zeta (\omega) 
  {\bf g}( \omega , {\bf i}(\omega), {\bf j}(\omega)).
$$
5) the global outcome is:
$$
\int_0^{+\infty} \gamma_t \, \bk (t) dt
$$
where $\gamma_t$ is the payoff at time $t$.\\
%Note that in ${\mathcal G}_\Pi$
%the players actually use mixed actions and the transition depends upon their realization.
%However on the discrete recursive  equation the expression is the same as if the game was played in expectation,  see Section \ref{Fleming}.
In  ${\mathcal G}_\Pi$ the state is deterministic and at each  time $t_n$ the players  know $\zeta_{t_n}$ and choose ${\bf i}_n$ (resp. ${\bf j}_n$). 

Consider now the mixed extension ${\mathcal G}^{II}_\Pi$ (Section 5) and let ${\mathcal V} _{\Pi}(t, \zeta)$  be the associated value.
% corresponding to the partition $\Pi$ and the evaluation $\bk$.

\begin{pro} $ $ \\
The value  ${\mathcal V} _{\Pi }(t, \zeta)$ satisfies the  recursive equation (\ref{RF1'}).
% \begin{eqnarray}\label{AX} 
%\hskip 2cm
%{\mathcal V} _{\Pi , \bk}(t_{n-1}, \zeta_{t_{n-1}}) &=& \val_{  {\bf X} \times{\bf Y }}  \,  [\int_{t_{n-1}}^{t_{n} } \langle \zeta_s , {\bf g}( ., {\bf x}(.),   {\bf y}(.))\rangle  \bk(s)  ds  + 
%{\mathcal V} _{\Pi , \bk}(t_{n}, \zeta_{t_{n}})
%\end{eqnarray}
\end{pro}
\noindent\underline{Proof}\\
The mixed action set for player 1 is $\tilde { \bf X}$ but due to the separability in $\omega$ one can work with $\bf X$.
Then it is easy to see that equation (\ref{RF3"}) corresponds to (\ref{RF1'}).
\square  $ $ \\

The analysis in section 5 thus implies that :\\
- any accumulation point  $U$ of the sequence ${\mathcal V} _{\Pi }$  is a viscosity solution of 
\begin{equation}\label{basic21}
\hskip 1,5cm 0 =  \frac{d}{dt} U (t , \zeta)  + \val_{  {\bf X} \times{\bf Y }} \; [\langle \zeta, {\bf g}( ., {\bf x}(.),   {\bf y}(.))\rangle  \bk(t) +  \langle  f (\zeta,  {\bf x}, {\bf y}), \nabla  U(t , \zeta)   \rangle ]
\end{equation}
- Equation (\ref{basic21}) has a unique viscosity solution.\\

%\begin{defi}
%A  real function $U$ on $ \rit^+ \times \Delta( \Omega ) $ is a viscosity solution of:
%\begin{equation}\label{basic21}
%\hskip 1,5cm 0 =  \frac{d}{dt} U (t , \zeta)  + \val_{  {\bf X} \times{\bf Y }} \; [\langle \zeta, {\bf g}( ., {\bf x}(.),   {\bf y}(.))\rangle  \bk(t) +  \langle  f (\zeta,  {\bf x}, {\bf y}), \nabla  U(t , \zeta)   \rangle ]
%\end{equation}
%if for any  real function $\phi$, ${\mathcal C}^1$ on $ \rit^+ \times \Delta( \Omega ) $ with $\phi -U$ having a strict minimum at $(\bar t,  \bar \zeta) \in \rit^+ \times \Delta( \Omega ) $
% \begin{equation}
% 0 \leq  \frac{d}{dt} \phi (\bar t, \bar \zeta)  + \val_{  {\bf X} \times{\bf Y }} \; [\langle \bar \zeta, {\bf g}( ., {\bf x}(.),   {\bf y}(.))\rangle  \bk(t) +  \langle  f (\bar \zeta,  {\bf x}, {\bf y}), \nabla \phi (\bar t , \bar \zeta)   \rangle ],
%\end{equation}
%and the dual condition.
%\end{defi}
%
%\begin{pro}$ $ \\ 
%Equation (\ref{basic21}) has a unique viscosity solution.
%\end{pro}
%\noindent\underline{Proof}\\
%Standard from differential games.
%\square\\

In particular let   $ \U  (t , \zeta)  = \langle \zeta, u(t, .) \rangle = \sum_{\omega} \zeta({\omega})  u(t, {\omega})$  where $u \in {\bf V} $.
%  is an accumulation point  of the family of values $\{v_{\Pi } \}$.

\begin{pro}$ $ \\ 
$ \U  (t , \zeta)$ is the  viscosity solution of (\ref{basic21}).
\end{pro}
\noindent\underline{Proof}\\
Follows from the fact that $ V_\Pi$  and ${\mathcal V} _{\Pi }$ satisfy the same recursive formula, hence $\U$ is an acumulation point  of the sequence ${\mathcal V} _{\Pi }$.
\square\\

This leads to the convergence property.

   \begin{cor}
   Both families $V_\Pi$ and $v_\Pi$ converge to some $V$ and $v$ with
   $$
   V(t, \zeta) = \sum_\omega \zeta(\omega) v(t, \omega).
   $$
   $V$ is the viscosity solution of (\ref{basic21}).\\
$v$    is the  viscosity solution of (\ref{ME1}).
%\begin{equation*}
%0 = \frac {d}{dt} v( t, z)  + \val_{X\times Y} \{{\bf g}(z, x, y) k(t) + {\bf q} ( x, y)[z, .] \circ v (t, \cdot)\}.
%\end{equation*}
   \end{cor}
   \noindent\underline{Proof}\\
   One has $\nabla  {\mathcal V} (t , \zeta)  = \{ v(t, .)\}$ hence:
\begin{eqnarray*}
0&=& \langle \zeta,  \frac{d}{dt}v (t ,.) \rangle + 
 \val_{  {\bf X} \times{\bf Y }} \; [\langle \zeta, {\bf g}( ., {\bf x}(.),   {\bf y}(.))\rangle  \bk(t) +  \sum_z [ \sum_\omega  \zeta (\omega) {\bf q}( {\bf x}(\omega), {\bf y}(\omega)) [ \omega, z] v (t , z)    ] \\
&=& \langle \zeta,  \frac{d}{dt}v (t ,.) \rangle + 
 \val_{  {\bf X} \times{\bf Y }} \; [\langle \zeta, {\bf g}( ., {\bf x}(.),   {\bf y}(.))\rangle  \bk(t) +  \sum_{\omega} \zeta (\omega)[ \sum_z  {\bf q}( {\bf x}(\omega), {\bf y}(\omega)) [ \omega, z]  v (t , z)    ]
\end{eqnarray*}
This gives: 
\begin{equation*}
0 = \langle \zeta,  \frac {d}{dt} v( t, .)  + \val_{X\times Y} [{\bf g}(., x, y) \bk(t) + {\bf q} ( x, y) [., \bullet]\circ v (t, \bullet)]\rangle
\end{equation*}
which is equivalent  to:
\begin{equation*}
0 = \frac {d}{dt} v( t, z)  + \val_{X\times Y} [{\bf g}(z, x, y) \bk(t) + {\bf q} ( x, y) [z,.] \circ v (t, \cdot)].
\end{equation*}
and this is (\ref{ME1}).

   \square

%%%%%%%%%%%%%%%%%%%%%%%%%%%%%%%%%%%%%%%%%%%%%%%%%%%%%%%%%%%%%%%%

\subsection{Stationary case}$ $ \\
We consider the case $\bk(t) = \rho e^{- \rho t}$ and again the  game along the  partition $\Pi$.

\subsubsection{Recursive formula}$ $ \\
The general recursive formula  (\ref{RF1})  takes now the following form:
\begin{pro}
 \begin{eqnarray}
 \nonumber  v_{\Pi,\rho}(t_{n}, Z_{t_{n}}) &=& \val_{X\times Y}  \,  \E_{z, x ,y} [\int_{t_{n}}^{t_{n+1} }{\bf g}(Z_s,  i, j)  \rho e^{- \rho }ds  + 
 v_{\Pi, \rho}(t_{n+1}, Z_{t_{n+1}} )]\\
 &=&
  \val_{X\times Y}  \, [ \E_{ z, x, y} (\int_{t_{n}}^{t_{n+1} }{\bf g}(Z_s,  i, j)  \rho e^{- \rho }ds ) +    \P^{\delta_n} ( x, y)[Z_{t_{n}}, .] \circ  v_{\Pi, \rho}(t_{n+1}, .)]
\end{eqnarray}
and if $\Pi$ is uniform, $v_{\Pi,\rho}(t, z) = e^{- \rho t} \nu_{\delta, \rho}(z)$ with:
\begin{equation}\label{RF1"} 
 \qquad \nu_{\delta, \rho}(Z_{0}) = \val_{X\times Y}  \, [\E_{z, x, y} (\int_{0}^{\delta }{\bf g}( Z_s, x, y)  \rho e^{- \rho }ds)  +  e^{-\rho \delta}\,   \P^{\delta}( x, y) [Z_0 ,.] \circ \nu_{\delta, \rho}(.)  ]
\end{equation}
\end{pro}

\subsubsection{Main equation}$ $ \\
The next result is standard, see e.g.  Neyman  \cite{N2013}, Prieto-Rumeau and Hernandez-Lerma \cite{PRHL}, p. 235. 
%check Guo Hernandez-Lerma \cite{GH1} \cite{GH2}.
We provide a short proof for convenience.
 \begin{pro} $ $ \\
  1)  For any $R \in {\mathcal M}$ and any $\rho \in (0, 1]$ the equation,  with variable $\varphi$ from $\Omega $ to $\rit$: 
\begin{equation}\label{ME1"}
 \rho \;  \varphi ( z) = \val_{X\times Y} [ \rho \; {\bf g}( z, x, y)   +  R( x, y)[z, .]  \circ  \; \varphi(. )] 
\end{equation}
has a  unique solution, denoted  $W_\rho$.\\
2) For any $\delta \in (0, 1]$ such that $\|  \delta R/  (1 -  \delta \rho)\| \leq 1$ the solution of (\ref{ME1"}) is the value of the repeated stochastic game with payoff $g$, transition $  P = I +  \delta R /  (1 -  \delta \rho)$ and discounted factor
$ \delta \rho$. 
 \end{pro}
\noindent\underline{Proof}\\
Recall from Shapley \cite{Sh}, that  the value $W_{\rho \delta} $ of a  repeated stochastic game  with payoff $\bf g$ and discounted factor $ \delta \rho$ satisfies:
 \begin{equation}\label{pasic}
W_{\rho \delta} (z) = \val_{X\times Y}[\delta \rho \, {\bf g}( z, x, y) + (1 -  \delta \rho) \E_{z, x, y} \{W_{\rho \delta} (.) \} ].
\end{equation}
Assume the  transition  to be of the form $ P = I + \delta q$ with $q \in \mathcal M$. One obtains:
\begin{equation}
W_{\rho \delta} (z) = \val_{X\times Y}[\delta \rho \, {\bf g}( z, x, y) + (1 -  \delta \rho)\{W_{\rho \delta} (z)  + \delta \,  q( x, y)[z,.] \circ  W_{\rho \delta} (.) \} ]
\end{equation}
which gives: 
\begin{equation}
\delta \rho\,  W_{\rho \delta} (z) = \val_{X\times Y}[\delta \rho Ý\, {\bf g}( z, x, y)   + \delta (1 -  \delta \rho)\,  q( x, y)[z,.] \circ  W_{\rho \delta} (.) ]
\end{equation}
so that: 
\begin{equation}
 \rho \, W_{\rho \delta} (z) = \val_{X\times Y}[\rho \, {\bf g}( z, x, y)   + (1 -  \delta \rho)\,  q( x, y)[z,.] \circ  W_{\rho \delta} (.) ].
\end{equation}
Hence with $q = R/  (1 -  \delta \rho)$  one obtains: 
\begin{equation}\label{music}
 \rho \,  W_{\rho \delta} (z) = \val_{X\times Y}[\rho \,  {\bf g} ( z, x, y)   + R ( x, y)[z,.] \circ  W_{\rho \delta} (.) ].
\end{equation}
\square

\subsubsection{Convergence}$ $ \\
Again the following result can be found in  Neyman \cite{N2013}, Theorem 1, see also Guo Hernandez-Lerma \cite{GH1, GH2}.
 
 \begin{pro}$ $ \\
As the mesh  $ \delta$ of the partition  $\Pi$ goes to 0,
$v_{\Pi,\rho}$ converges to the solution  $W_{\rho}$ of  (\ref{ME1"}) with $R = {\bf q}$:
 \begin{equation}\label{basic11}
 \rho \; W_{\rho}( z) = \val_{X\times Y} [ \rho \; g( z, x, y)   +    {\bf q}( x, y)[z, .] \circ  \;  W_{\rho} (. )] 
\end{equation}

 \end{pro} 
 \noindent\underline{Proof}\\
 Consider the strategy $\sigma$  of Player 1  in $G_{\Pi}$ defined as follows:  at state $z$,   use an optimal strategy  $x \in X = \Delta(I)$,  for $W_{\rho}(z)$ given by   (\ref{basic11}).
 Let us evaluate, given $\tau$, strategy of Player 2,  the following amount: 
 $$
 A_1 =   \E_{\sigma, \tau}[\int_{t_1}^{t_{2} }g_\Pi (s)  \rho e^{ -\rho s} ds  +  e^{ -\rho t_1}  W_{\rho}( Z _{t_1}) ].
 $$
Let $x_1$ the mixed move of Player 1 at stage one given $Z_0= {\hat Z}_1$. Then if $y_1$ is induced by $\tau$, there exists a constant $L$ such that:
\begin{eqnarray*}
A_1 &\geq& \ \delta_1 \rho \; g( {\hat Z}_1, x_1, y_1) + ( 1 - \delta_1 \rho ) [ W_{\rho}({\hat Z}_1) + \delta_1 \; {\bf q} ( x_1, y_1 )[{\hat Z}_1,.] \circ W_{\rho}(.)]  - \delta_1  L   \delta \\
&\geq& \ \delta_1 \rho \; g( {\hat Z}_1, x_1, y_1)  - \delta_1 \rho \;  W_{\rho}({\hat Z}_1) +  \delta_1 \; {\bf q} ( x_1, y_1 )[{\hat Z}_1,.]  \circ W_{\rho}(.)   + W_{\rho}({\hat Z}_1) - 2  \delta_1L   \delta\\
&\geq&W_{\rho}({\hat Z}_1) - 2 \delta_1 L   \delta.
\end{eqnarray*}
Similarly  let: 
 $$ 
 A_{n}=  \E_{\sigma, \tau}[\int_{t_{n}}^{t_{n+1} }g_\Pi (s)  \rho e^{ -\rho s} ds  ds  + e^{- \rho t _n} W_{\rho}( {\hat Z}_{n+1}) | h_n]
 $$
 where $h_n = ( {\hat Z}_1, i_1, j_1, \cdots, i_{n-1}, j_{n-1}, {\hat Z}_n)$.\\
 Then, with obvious notations: 
\begin{eqnarray*}
A_n &\geq& 
e^{- \rho t _{n-1}}   [ \delta_n \rho \; g( {\hat Z}_n, x_n, y_n) + ( 1 - \delta_n \rho ) [ W_{\rho}({\hat Z}_n) + \delta_n\; {\bf q} ( x_n, y_n )[{\hat Z}_n,.]  \circ W_{\rho}(.) - \delta_n L   \delta  ] \\
&\geq& e^{- \rho t _{n-1}}  [ \delta_n \rho \; g( {\hat Z}_n, x_n, y_n)  - \delta_n \rho\;  W_{\rho}({\hat Z}_n) +  \delta_n \; {\bf q} ( x_n, y_n )[{\hat Z}_n,.]  \circ W_{\rho}(.)   + W_{\rho}({\hat Z}_n) -2  \delta_n  L   \delta ]\\
&\geq& e^{- \rho t _{n-1}} [W_{\rho}({\hat Z}_n) - 2 \delta_n  L  \delta ].
\end{eqnarray*}
 Taking the sum and the expectation, one obtains that the payoff induced by $(\sigma, \tau)$ in $G_\Pi$ satisfies:
 \begin{equation*}
\E_{\sigma, \tau}[\int_{0}^{+\infty} g_\Pi (s)  k(s) ds]  \geq W_{\rho} ({\hat Z}_1) - 2( \sum_n  \delta_n e^{- \rho t _{n-1}} ) L  \delta
\end{equation*}
and $ ( \sum_n  \delta_n e^{- \rho t _{n-1}} )  L  \delta  \rightarrow 0$ as $  \delta  \rightarrow 0$.
\square\\

 \noindent{\bf Comments}:
 
  The proof in Neyman \cite{N2013} is done, in the finite case,  for a uniform partition but shows the robustness with respect to the parameters (converging family of games). 
  
This procedure of proof is reminiscent of the ``direct approach" introduced by Isaacs \cite{I}. To show convergence of the  family of values of the discretizations $v_\Pi$: $ i)$  one identifies a tentative limit value $v$  and a recursive formula $RF(v)$ and $ ii)$ one shows that to  play in the discretized  game $G_\Pi$ an optimal strategy in $RF(v)$  gives an amount close to $v$ for $\delta$ small enough.

  For an alternative approach  and proof, based on properties of  the Shapley operator, see Sorin and Vigeral \cite{SV3}.
  
 Remark that if $\bk(t) = \rho e^{-\rho t}$,  $v ( t,z) = e^{-\rho t} \nu (z)$ satisfies (\ref{ME1})  iff $\nu(z)$ satisfies (\ref{basic11}).

 %%%%%%%%%%%%%%%%%%%%%%%%%%%%%%%%%%%%%%%%%
  %%%%%%%%%%%%%%%%%%%%%%%%%%%%%%%%%%%%%%%%%
   %%%%%%%%%%%%%%%%%%%%%%%%%%%%%%%%%%%%%%%%%

%%%%%%%%%%%%%%%%%%%%%%%%%%%%%%%%%%%%%%%%%%%%%%%%%%%%%%%%%%%%%%%%%%%%%%%%%%%%%%%%%%%%%%%%%%%%%%%%%%%%%%%%%%%%%%%%%%%%%%%%%%%%%%%%%%%%%%BBBBBBBBBBBBBBBBBBBBBBBBBBBBBBBBBBBBBBBBBBBBBBBBBBBBBBBBBBBB%%BBBBBBBBBBBBBBBBBBBBBBBBBBBBBBBBBBBBBBBBBBBBBBBBBBBBBBBBBBBB
%%BBBBBBBBBBBBBBBBBBBBBBBBBBBBBBBBBBBBBBBBBBBBBBBBBBBBBBBBBBBB

\section{State controlled and not observed}

This section studies the game $\G $  where  the  process $Z_t$ is controlled by both players but not observed. However  the past actions are known:  this defines a symmetric framework were  the new state  variable is the law of $Z_t$,  $\zeta_t \in \Delta(\Omega)$. 
%Denote by  $\G (\zeta)$ the corresponding game and let $\V(\zeta)$ be its value.
Even in the stationary case  there is no explicit  smooth solution to the main equation  hence a direct approach  for proving   convergence, as in the previous Section 3.2,  is not feasible.  
 
 Here also the analysis will be trough the connection to  a differential game  $\overline{\mathcal G}$ on $\Delta(\Omega)$ but different from the previous one $\mathcal G$, introduced in Section 3.
 
Given a partition $\Pi$ denote by $\G_\Pi$ the associated game and again, since $\bk$ is fixed during the analysis we will write $\V_{\Pi}$ for  its value $\V_{\Pi, \bk}$ defined on $ \rit ^+\times\Delta (\Omega)$.\\
Recall that given  the initial law $ \zeta_{t_{n}}$ and  the actions $(i_{t_{n}}, j_{t_{n}})= ( i,  j)$  one has:
\begin{equation}\label{hic}
 \zeta_{t_{n+1}}^{ij} = \zeta _{{t_{n}}+\delta_n}= \zeta_{t_{n}} \ast  \P^{\delta_n } ( i,  j)
\end{equation}
  and that  this parameter is known by both players.\\
  Extend ${{\bf g} (., x, y)}$  from  $\Omega$ to  $\Delta (\Omega)$   by linearity:   ${\bf g} (\zeta, x, y)= \sum \zeta(z) {\bf g} (z, x, y)$.\\
 %Again, since $\bk$ is fixed during the analysis we will write $\V_{\Pi}$ for $\V_{\Pi, \bk},$ defined on $ \rit ^+\times\Delta (\Omega)$.

 \subsection{Recursive formula}$ $ \\
In this framework the recursive structure leads to:
 \begin{pro} $ $ \\
The value  $\V_{\Pi}$ satisfies the following recursive formula:
 \begin{eqnarray}\label{RF2} 
\V_{\Pi}(t_{n}, \zeta_{t_{n}}) &=& \val_{X\times Y}  \,  \E_{\zeta, x, y} [\int_{t_{n}}^{t_{n+1} }{\bf g}(\zeta_s,  i, j)  \bk(s) ds  + 
\V_{\Pi}(t_{n+1}, \zeta_{t_{n+1}}^{ij} )]
 %&=& \val_{X\times Y}  \, [\int_{t_{n-1}}^{t_{n} }{\bf g}( \zeta_s,  x, y) k(s) ds  +   \E_{x,y} v_{\Pi}(t_{n},   \P^{\delta_n} (\zeta_{t_{n-1}}, i, j)  )]
\end{eqnarray}

\end{pro}
\noindent\underline{Proof}\\
Standard, since $\G_\Pi$ is basically a stochastic game with parameter $\zeta$.
\square\\

\subsection{Main equation}$ $ \\
Consider the differential game $\overline { \mathcal G}$ on $\Delta(\Omega)$ with actions sets $I$ and $J$,   dynamics on $\Delta (\Omega) \times \rit^+$ given by: 
$$
\dot \zeta_t = \zeta_t  * {\bf q} ( i, j), 
$$
 current payoff   ${\bf g} (\zeta, i,  j)$ and evaluation $\bk$.

As in Section 5, consider the discretized mixed extension $\overline { \mathcal G}^{II}_\Pi$  to $X\times Y$ and let $\overline{\mathcal V}_\Pi$ be its value.
  \begin{pro}$ $ \\ 
$\overline{\mathcal V}_\Pi$ satisfies (\ref{RF2}).
\end{pro}
\noindent\underline{Proof}\\
$\overline{\mathcal V}_\Pi$ satisfies (\ref{RF3"}) which is, using (\ref{hic}),  equivalent  to (\ref{RF2}).
\square\\

  The analysis in Section 5, Proposition  \ref{fin}  thus implies:
\begin{pro}$ $ \\ 
The family of values  $\V_{\Pi}$ converge to $\V$ unique viscosity solution of : 
\begin{equation}\label{MF2}
0 = \frac {d}{dt} u( t,  \zeta)  + \val_{X\times Y} [{\bf g}(\zeta, x, y) \bk(t) + \langle \zeta * {\bf q} ( x, y), \nabla  u (t, \zeta)].
\end{equation}

\end{pro}\subsection{Stationary case}$ $ \\
Assume $\bk(t) = \rho e^{- \rho t}$.\\% and the  game along the  partition $\Pi$.
In this case one has $ \V( \zeta, t) = e^{-\rho t} \v (\zeta)$ hence (\ref{MF2}) becomes
\begin{equation}\label{NR}
 \rho \v ( \zeta)   =  \val_{X\times Y} [ \rho \;  {\bf g}(\zeta, x, y) + \langle \zeta * {\bf q} ( x, y), \nabla  \v ( \zeta)\rangle]
\end{equation}

\subsection{Comments}$ $ \\
A  differential game similar to $\overline { \mathcal G}$ where the state space is the set of probabilities on some set $\Omega$ has been studied in full generality by Cardaliaguet and Quincampoix \cite{CQ}, see also As Soulaimani  \cite{AS2008}.\\
 Equation (\ref{NR})  is satisfied by the value of the Non-Revealing game in the framework analyzed by  Cardaliaguet, Rainer, Rosenberg and Vieille \cite{CRRV} see Section 6.
 
 %%%%%%%%%%%%%%%%%%%%%%%%%%%%%%%%%%%%%%%%%%%%%%%%%%%%%%%%%%%%%%%%%%%%%%%%%%%SECTION 5 %%%%%%%%%%%%%%%%%%%%%%%%%%%%%%%%%%%%%%%%%%%%%%%%%%%%%%%%%%%%%%%%%%%%%%%%%%%%%%%%%%%%%%%%%

\section{Discretization and mixed extension of differential games}

  We study here the value of a continuous time game by introducing a   time discretization 
$\Pi$ and analyzing the limit behavior  of the  associated family of values $v_\Pi$ as the mesh of the partition vanishes. This approach was initiated  in Fleming \cite{F57},  \cite{F61},\cite{F64}, and developped in Friedman \cite{Fri71}, \cite{Fri74},  Eliott and Kalton \cite{EK}.

A differential game  $\gamma$ is defined trough the following components:
$Z \subset \rit^n$ is the state space, $I$ and $J$ are the action sets of player 1 (maximizer) and 2, $f $ from $Z\times  I\times J$ to $\rit^n$ is  the dynamics kernel,
$g$  from $Z\times  I\times J$ to  $\rit$ is the payoff-flow  function and $k$ from $\rit ^+$ to $\rit^+$ determines the evaluation.\\
Formally the  dynamics is defined on $[0, + \infty)\times Z$ by  :
\begin{equation}\label{Dy}
\dot z_t = f(z_t, i_t, j_t) 
\end{equation}
 and the total payoff is :
 $$
 \int_0^{+\infty}g( z_s, i_s, j_s) \, k(s) \, ds.
 $$
 We assume:\\
  $I$ and $J$  metric compact sets,\\
 $f$ and $g$ continuous and uniformly Lipschitz in $z$,\\
  $g$ bounded,\\
 $k$ %: [0, +\infty)  \rightarrow  [0, +\infty)  
   Lipschitz with $\int_0^{+\infty} k(s) ds = 1$.\\
$ \Phi^h(z; i,j ) $ denote  the value at time $t+h$ of the solution of (\ref{Dy}) starting at time $t$ from $z$ and with  play $\{ i_s = i, j_s = j\}$ on $[t, t+h]$.\\
To define the strategies we have to specify the information:  we assume  that the players know  the initial state $z_0$, and at time $t$ the previous play
 $\{ i_s, j_s; 0\leq s < t\}$ hence the trajectory of the state $\{ z_s; 0\leq s\leq t\}$.\\
 The analysis below will show that Markov strategies (i.e. depending only, at time $t$, on $t$ and $z_t$)  will suffice.

\subsection{Deterministic analysis}$ $ 

Let $\Pi = (\{t_n\} , n= 1, ...)$  be a partition of  $[0, + \infty)$  with $t_1 = 0,  \delta_n = t_{n+1} - t_n$ and $ \delta = \sup \delta_n$. We consider the associated discrete time game  $\gamma_\Pi$ where on each interval  $[t_n, t_{n+1})$ players use constant moves $(i_n,j_n)$  in  $I\times J$.  This defines the dynamics on the state. At time $t_{n+1}$,  $(i_n,j_n)$ is announced and the  corresponding  value of the state, $z_{t_{n+1}}  = \Phi^{\delta_n}(z_{t_{n}};  i_n, j_n) $ is known. \\
The associated maxmin  $w_\Pi ^-$ satisfies the recursive formula: 
\begin{equation}\label{RF3}
w^-_\Pi(t_{n}, z_{t_{n}} )= \sup_I \inf_J   [\int_{t_{n}}^{t_{n+1} }{g}(z_s,  i, j)  k(s) ds  +  
 w^-_{\Pi}(t_{n +1}, z_{t_{n+1}} )]\\
\end{equation}
The fonction $w^-_\Pi(.,z) $
 is extended by linearity to $[0, + \infty)$ and note that: 
\begin{equation}\label{unif}
\forall \varepsilon>  0, \exists T, \mbox{ such that $ t \geq T $ implies }   |w^-_\Pi (t, .) |\leq \varepsilon
\end{equation}
and that all ``value'' functions  that we will consider here will satisfy this property. \\

The next four results follow from  the analysis in Evans and Souganidis  \cite{ES}, see also Barron, Evans and Jensen  \cite{BEJ}, Souganidis  \cite{ Sou85}  and the presentation in Bardi and Capuzzo-Dolcetta \cite{BC}, Chapter VII, Section 3.2.

\begin{pro}$ $ \\ 
The family $\{ w^-_\Pi(t,z) \}$ is  uniformly equicontinuous in both variables. 
\end{pro}
Hence the set $U$ of accumulation points of the family  $\{w_\Pi ^-\}$ (for the uniform convergence on compact subsets of $\rit^+ \times Z$),  as  the mesh $\delta$ of $\Pi$ goes to zero, is non empty.\\

We first introduce  the notion of viscosity solution, see Crandall and Lions  \cite{CL}.

\begin{defi} Given an Hamiltonian $H$ from $\rit^+\times Z\times \rit^n$ to $\rit$, a  continuous real function $u$ on  $\rit ^+\times Z $ is a viscosity solution of:
\begin{equation}\label{ME}
0 = \frac {d}{dt} u( t, z)  + H(t, z, \nabla u (t, z) )
\end{equation}
if for any  real function $\psi$, ${\mathcal C}^1$ on $\rit ^+\times Z $ with $u - \psi$ having a strict maximum at $(\bar t,  \bar z)\in \rit ^+\times Z $:
 \begin{equation*}
0 \leq  \frac {d}{dt} \psi( \bar t, \bar z)  +  H(t, z, \nabla \psi (t, z) )
\end{equation*}
and the dual condition holds.
\end{defi}

We can now introduce the Hamilton-Jacobi-Isaacs (HJI) equation that follows from (\ref{RF3}),  corresponding to the Hamiltonian:
\begin{equation}\label{h^-}
h^-( t, z, p) = \sup_I \inf_J  [{g}(z, i, j ) k(t)   +   \langle f(z, i, j ), p  \rangle ].
\end{equation}

\begin{pro}$ $ \\ 
Any accumulation point $u \in U$  is a viscosity solution of: 
\begin{equation}\label{ME3}
0 = \frac {d}{dt} u ( t, z)  +  \sup_I \inf_J  [{g}(z, i, j ) k(t)   +   \langle f(z, i, j ), \nabla u (t, z) \rangle ].
\end{equation}
\end{pro}

Note that in the discounted case,   $k(t) = \lambda e^{- \lambda t}$,   with the change of variable $u(t, z) =  e^{- \lambda t} \phi (z) $, one obtains:
\begin{equation}
\lambda \phi (z) = \sup_I \inf_J  [ \lambda {g}(z, i, j ) +   \langle f(z, i, j ), \nabla \phi(z) \rangle ].
\end{equation}
The main property is the following:
\begin{pro}$ $ \\ 
Equation (\ref{ME3}) has a unique viscosity solution.
\end{pro}
Recall that this notion and this result are  due to Crandall and Lions  \cite{CL}, for more  properties see Crandall, Ishii and Lions \cite{CIL}.\\

The uniqueness of accumulation point implies:
\begin{cor}$ $ \\ 
 The family $\{w_\Pi ^-\}$ converges to some $w^-$.
\end{cor}

An alternative approach is the consider the game $\gamma$  in normal form on $\rit^+$. Let  $w_\infty ^- $ be the  maxmin (lower value) of the continuous time differential game   played using  non anticipative strategies with delay. 
Then from Evans and Souganidis  \cite{ES},  extended in  Cardaliaguet \cite{Cii10}, Chapter 3,  one obtains:
\begin{pro}$ $ \\ 
1) $w_\infty ^-$  is a viscosity solution of  
(\ref{ME3}).\\
2) Hence:
$$
w_\infty ^-  = w^-.
$$
\end{pro}

Obviously  similar properties hold for  the minmax $w_\Pi ^+$ and $ w_\infty ^+$. \\

Finally define  Isaacs's condition on $I \times J$ by : 
\begin{equation} \label{Is1}
 \sup_I \inf_J  [{g}(z, i, j )\;  k(t) +   \langle f(z, i, j ), p \rangle ] =
   \inf_J \sup_I [{g}(z, i, j )\;  k(t) +   \langle f(z, i, j ), p \rangle ],   \; 
  \forall  t \in \rit^+,  \forall z \in Z, \forall p \in \rit^n,
\end{equation}
which, with the notation (\ref{h^-}), corresponds to :
$$
h^-(t, z, p) = h^+ (t, z, p).
$$

 \begin{pro}$ $ \\
 Assume condition (\ref{Is1}).\\
 Then the limit value exists, in the sense that:
 $$
 w^- = w^+ ( = w_\infty^- = w_\infty^+) 
 $$
\end{pro}

Note that the same analysis holds if the players use strategies that depend only at time  $t_n$ on $t_n$ and $z_{t_n}$.
 %%%%%%%%%%%%%%%%%%%%%%%%%%%%%%%%%%%%%%%%%%%%%%%%%%%%%%%%%%%%%%%%%%%%%%%%%%%%%%%%%%%%%%%%%%%%%%%%%%%%%%%%%%%%%%%%%%%%%%%%%%%%%

 \subsection{Mixed extension}$ $ 
 
 We define two mixed extensions of $\gamma$ as follows:
for each partition $\Pi$ we introduce  two  discrete time games  associated to $\gamma_\Pi$  and played  on $X = \Delta (I)$ and $Y = \Delta (Y)$  (set  of probabilities on $I$ and $J$ respectively). We will then prove that their asymptotic properties 
coincide.

 \subsubsection {Deterministic actions}$ $
  \\
  The first game $\Gamma^I$
  is defined as in subsection 5.1 were $X$ and $Y $ are now the sets of actions (this corresponds to ``relaxed controls") replacing $I$ and $J$. \\
The  main point is that the dynamics  $f$  (hence the  flow) and the payoff  $g$ are extended   to $X\times Y$ by taking  the expectation w.r.t. $x$ and $y$:
$$
f(z,x,y) = \int_{I\times J} f(z, i, j) x(di) y(dj)
$$
\begin{equation}\label{Dy'}
\dot z_t = f(z_t, x_t, y_t) 
\end{equation}
$$
g(z,x,y) = \int_{I\times J} g(z, i, j) x(di) y(dj).
$$
$\Gamma_\Pi^I$ is the associated discrete time game   where on each interval  $[t_n, t_{n+1})$ players use constant actions $(x_n,y_n)$  in  $X\times Y$.  This defines the dynamics: $\bar \Phi^h(z; x, y) $ denotes  the value at time $t+h$ of the solution of (\ref{Dy'}) starting at time $t$ from $z$ and with  play $\{ x_s = x, y_s = y\}$ on $[t, t+h]$. Note that  $ \bar \Phi^h(z; x, y) $ is {\it not} the bilinear extension of  $ \Phi^h(z; i,j ) $. At time $t_{n+1}$,   $(x_n,y_n)$ is announced and the current  value of the state, $z_{t_{n+1}}  = \bar\Phi^{\delta_n}(z_{t_{n}};  x_n, y_n) $ is known. \\
The maxmin $W_\Pi ^-$ satisfies the recursive formula:
$$
W^-_\Pi(t_{n}, z_{t_{n}} )= \sup_X \inf_Y  [\int_{t_{n}}^{t_{n+1} }{g}(z_s,  x, y)  k(s) ds  +  
 W^-_{\Pi}(t_{n +1}, z_{t_{n+1}} )].
$$
 The analysis of the previous paragraph applies,  leading to:

\begin{pro}$ $ \\ 
The family $\{W^-_\Pi(t,z)  \}$ is uniformly equicontinuous in both variables. 
\end{pro}
The HJI equation corresponds here to  the Hamiltonian:
\begin{equation}\label{H^-}
H^-( t, z, p) = \sup_X \inf_Y \,  [{g}(z, x, y ) k(t)   +   \langle f(z, x, y ), p \rangle ].
\end{equation}

\begin{pro}$ $ \\ 
1) Any accumulation point of the family  $\{W_\Pi ^-\}$, as  the mesh $\delta$ of $\Pi$ goes to zero, is a viscosity solution of: 
\begin{equation}\label{ME3-}
0 = \frac {d}{dt} W^-( t, z)  +  \sup_X\inf_Y  [{g}(z, x, y ) k(t) +   \langle f(z, x, y ), \nabla W^- (t, z) \rangle ]
\end{equation}
2) The family $\{W_\Pi ^-\}$ converges to   $W^-$, unique viscosity solution of (\ref{ME3-}).
\end{pro}

Finally let  $W_\infty ^- $ be the  maxmin  of the differential game  $\Gamma^I$ played  (on $X\times Y$) using  non anticipative strategies with delay. 
Then:
\begin{pro}$ $ \\ 
1) $W_\infty ^-$  is a viscosity solution of  
(\ref{ME3-}).\\
2) 
$$
W_\infty ^-  = W^-.
$$
\end{pro}
As above, similar properties hold for $W_\Pi ^+$ and $ W_\infty ^+$. \\

Due to the bilinear extension,   Isaacs's condition on $X \times Y$ which is, with the notation (\ref{H^-}): 
\begin{eqnarray}\label{Is2}
H^{-}(t, z, p) = H^+ (t, z, p) 
% & & \sup_X \inf_Y  [{g}(z, x, y ) k(t) +   \langle f(z, x, y), p \rangle ] = \cr
% \nonumber &&   \inf_Y \sup_X  [{g}(z, x, y ) k(t) +   \langle f(z, x, y), p \rangle ],
 \quad \forall  t \in \rit^+,  \forall z \in Z, \forall p \in \rit^n,
\end{eqnarray}
always holds. Thus one obtains:

 \begin{pro}\label{R}$ $ \\
The limit value $W$ exists:
 $$
 W =W^- = W^+,
 $$
   and  is also the  value of the differential game played on $X \times Y$.\\
   It is the unique viscosity solution of :
%\end{pro}
  %Remark that due to  property (\ref{Is2}), (\ref{Ap2}) can be written as 
\begin{equation}\label{ME3'}
0 = \frac {d}{dt} W ( t, z)  + \val_{X\times Y}  [{g}(z, x, y) k(t) +   \langle f(z, x, y), \nabla W (t, z)\rangle ].
\end{equation}
\end{pro}$ $

%%%%%%%%%%%%%%%%%%%%%%%%%%%%%%%%%%

 \subsubsection {Random actions}$ $
  \\
  We define now another game $ {\Gamma}_\Pi^{II}$ where  the actions $(i_n,j_n) \in I\times J$  are chosen at random at time $t_n$ according to $x_n \in X$ and $y_n \in Y $,  then constant on $[t_n, t_{n+1})$   and announced at time $t_{n+1}$.  The new state is thus, if $(i_n, j_n) = (i,j)$, $z_{t_{n+1}}^{ij}  = \Phi^{\delta_n}(z_{t_{n}}; i, j ) $. \\
  It is clear, see e.g.  \cite{MSZ} Chapter 4,  that the next dynamic programming property holds: 
    \begin{pro}
The  game  $ {\Gamma}_\Pi^{II}$ has a value $\W_\Pi$, which satisfies the recursive formula:
\begin{equation}\label{RF3"}
\W_\Pi(t_{n}, z_{t_{n}} )= \val_{X\times Y} \E_{x,y}  [\int_{t_{n}}^{t_{n+1} }{g}(z_s,  i, j)  k(s) ds  +  
 \W_{\Pi}(t_{n +1}, z_{t_{n+1}}^{ij} )]
\end{equation}
\end{pro} 
\noindent and  given the hypothesis one obtains as above:
 \begin{pro}$ $ \\ 
The family $\{\W_\Pi(t,z)  \}$ is equicontinuous in both variables. 
\end{pro}
 Moreover one has:
\begin{pro}\label{fin}$ $ \\ 
1) Any accumulation point $ \U$ of the family  $\{\W_\Pi \}$, as  the mesh $\delta$ of $\Pi$ goes to zero, is a viscosity solution of (\ref{ME3'}).\\
2) The family $\{\W_\Pi \}$ converges to    $\W$, unique solution of (\ref{ME3'}).
\end{pro}
\noindent\underline{Proof}\\
1) Let $\psi ( t, z)$ be a  ${\mathcal C}^1$ test function such that  $\U -  \psi$ has a strict maximum at $(\bar t,  \bar  z)$. Consider a sequence $W_m = \W_{ \Pi(m)}$ converging uniformy locally  to $/U$  as $m \rightarrow \infty$ and  let $( t^*(m), z(m))$ be a minimizing sequence for  $  ( \psi - W_m) ( t,  z), t \in \Pi(m) $. In particular $( t^*(m), z(m)) $ converges to $(\bar t,  \bar z)$ as $m \rightarrow \infty$ . 
Given $x^*(m)$ optimal  in (\ref{RF3"}) one has  with  $t^*(m) = t_{n} \in \Pi(m)$: 
$$ 
W_m(t_{n}, z(m)) 
\leq 
\E_{x^*(m), y}  [ \int_{t_{n}}^{t_{n+1} }{ g}( z_s,  i, j) k(s) ds   +    W_m(t_{n+1}, z^{ij}_{t_{n+1}} )], \qquad \forall y \in Y 
 $$
 so that by the choice of $( t^*(m), z(m))$:
 \begin{eqnarray*}
\psi (t_{n}, z(m)) 
&\leq &
\E_{z,x^*(m), y}    [ \int_{t_{n}}^{t_{n+1} }{g}( z_s,  i, j) k(s) ds   +      \psi (t_{n+1}, z^{ij}_{t_n+1})]\\
%&\leq &
%\E_{x^*(m), y} [ \int_{t_{n}}^{t_{n+1} }{g}( z_s,  i, j) k(s) ds]   +  \sum_{i,j} x^{*i}(m) \, y^j \,     \psi (t_{n+1}, z_{t_{n}}  \P^{\delta_n } ( i, j))\\
  &\leq&
 \delta_n k(t_{n}) \; {g} (z(m), x^*(m), y) + \psi (t_{n+1}, z(m) ) \\
 &&+\,  \delta_n \;  \E_{ x^{*}(m) \, y} \,  \langle f( z(m),  i, j),   \nabla \psi (t_{n+1}, z(m))  \rangle + o(\delta_n).
\end{eqnarray*}
This implies:
\begin{eqnarray*}
0 &\leq &  \delta_n \;  \frac{d}{dt} \psi (t_{n}, z(m) )  +  \delta_n \; k(t_{n}) \; {g} (z(m), x^*(m), y) + \delta_n \;   \E_{ x^{*}(m) \, y}   \langle  f( z(m),  i, j) , \nabla \psi (t_n, z(m)) \rangle   + o(\delta_n)
\end{eqnarray*}
hence dividing  by $\delta_n$ and taking the limit  as $m \rightarrow \infty$ one obtains, for some accumulation point $x^* \in \Delta (I)$: 
\begin{equation}
0 \leq   \frac{d}{dt} \psi (\bar t, \bar z ) +  k(\bar t) \; {g} (\bar z, x^*, y)  +   \E_{ x^{*}\, y}  \langle f( \bar z,  i, j)  ,   \nabla \psi (\bar t, \bar z ) \rangle ,   \qquad \forall y \in Y.
\end{equation}
Thus $\U$ is  a viscosity solution of :
$$
0 = \frac {d}{dt} u( t, z)  + \val_{X\times Y}  \int_{I\times J} [{g}(z, i, j) k(t) +   \langle f(z, i, j ), \nabla u (t, z)\rangle ]x(di) y (dj)
$$
which by linearity, reduces to (\ref{ME3'}).\\
2) The proof of uniqueness follows from Proposition \ref{R}.
\square \\

Note again that the same analysis holds if the players use strategies that depend only at time  $t_n$ on $t_n$ and  $z_{t_n}$.

 \subsubsection {Comments}$ $
 
 Both games lead to the same limit PDE (\ref{ME3'})
 but with different  sequences of approximations:\\
 % In the case of incomplete information even the state variables appearing the programming dynamic formula  differ.
 In the first case ($\Gamma^I$),  the evolution is deterministic and the state (or $(x,y)$) is announced. \\
 In the second case  ($\Gamma^{II}$), the evolution is random and the state (or the actions)  are announced (the knowledge of $(x,y)$ would  not be enough).
 
 The fact that both games have same limit value is a justification for playing distribution or mixed actions as pure actions in continuous time  and  for assuming that the distributions   are observed, see Neyman \cite{N2012}.
 
%***\\
%
%another variant would have the dynamics defined by 
%$I\times J$ but only the expectation known. Then 
%$$
% z_h= \bar \Phi ( z, x, y ; h) = E_{x,y} [\Phi( z, i, j; h)
% $$
% 
% ****

Remark also  that  the same analysis holds if $f$  and $g$ depend in addition  continuously on $t$. 

A related study of differential games with mixed actions, but concerned with the analysis trough strategies can be found in \cite{BLQ}, \cite{BQRX}, \cite{JQX}.

The advantage of working with discretization is to have a well defined and simple  set of strategies hence the recursive formula is immediate to check for  the associated maxmin or minmax $W^{\pm}_\Pi$. 
On the other hand the main equation (HJI) is satisfied by accumulation points.\\
The  use of  mixed actions in extensions of type II allows to have values in the associated game.

\section{Concluding comments and extensions}

This research is part of an analysis of asymptotic properties  of dynamic games through their recursive structure : operator approach \cite{Sh}, \cite{RS}.

Recall that the analysis in terms of repeated games may lead to non convergence,  in the framework of Section 3 with compact action spaces, see Vigeral \cite{V4}, or in the framework of Section 4 even with finite action spaces, see Ziliotto \cite{Zi1} (for an overview of similar phenomena see Sorin and Vigeral \cite{SV2}).

The approach in terms of vanishing duration of a continuous time process allows, via 
the extension of the state space from $\Omega$ to $\Delta (\Omega)$   to  obtain smooth transition and nice limit behavior as $\delta$ vanishes. 

A similar procedure has been analyzed by Neyman \cite{N2013}, in the finite case, for more general classes of approximating games  and developed in Sorin and Vigeral \cite{SV3}.

The case of  private information on the state variable has been treated by Cardaliaguet, Rainer, Rosenberg and Vieille  \cite{CRRV} in the stationary finite framework: the viscosity solution corresponding to (ME) involves a geometric aspect due to the revelation of information that makes the analysis much more difficult. The (ME) obtained  here in  Section 3 corresponds to the Non Revealing value that players can obtained without using their private information.  
 
Let us  finally mention three directions of research:\\
the study of the general symmetric case i.e. a  framework between Section 3 and Section 4 where the players receive partially revealing symmetric signals on the state, \cite{S2016},\\
the asymptotic properties when both the evaluation tends to $+ \infty$ and the mesh goes to 0: in the stationary case this means both $\rho$ and $\delta$ vanishes. In the framework of Section 3, with finite actions spaces this was done by Neyman \cite{N2013} using the algebraic property of equation (\ref{basic11}),\\
the construction of optimal strategies based  at time $t$ on the current state $z_t$ and the instantaneous discount rate $\bk (t)/\sum_t^{+\infty} \bk (s) ds $.

\end{document}